\documentclass[11pt, a4paper]{article}
\usepackage{mathrsfs, amsmath, amsthm, amssymb, bm, mathtools, thm-restate}
\usepackage{empheq}
\usepackage{graphicx, xcolor, tikz}
\usetikzlibrary{arrows.meta, intersections, calc, shapes, snakes}
\usepackage{caption}
\usepackage[labelformat=simple]{subcaption}

\usepackage{array}
\usepackage{cases}
\makeatletter
\def\th@plain{%
  \upshape 
}
\makeatother

\makeatletter
\renewenvironment{proof}[1][\proofname]{\par
  \pushQED{\qed}%
  \normalfont \topsep6\p@\@plus6\p@\relax
  \trivlist
  \item[\hskip\labelsep
        \bfseries
    #1\@addpunct{.}]\ignorespaces
}{%
  \popQED\endtrivlist\@endpefalse
}
\makeatother

\newtheorem{theorem}{Theorem}[section]
\newtheorem{lemma}[theorem]{Lemma}
\newtheorem{corollary}[theorem]{Corollary}

\newtheorem*{conjecture*}{Conjecture}
\theoremstyle{definition}
\newtheorem{definition}{Definition}

\newtheorem{proposition}{Proposition}

\usepackage[left=25mm,top=25mm,bottom=25mm,right=25mm]{geometry}
\usepackage[T1]{fontenc}
\usepackage[inline]{enumitem}
\usepackage[square, numbers, sort&compress]{natbib}

\usepackage[
  bookmarks=true,%
  bookmarksnumbered=true, 
  bookmarksopen=true, 
  plainpages=false,%
  colorlinks=true, 
  linkcolor=blue, 
  citecolor=blue,%
  anchorcolor=green,
  urlcolor=blue,
  hyperindex=true]{hyperref}

\newcommand{\etal}{et~al.\ }

\begin{document}
\title{Weak degeneracy of planar graphs without 4- and 6-cycles}
\author{Tao Wang\thanks{Center for Applied Mathematics, Henan University, Kaifeng, 475004, China. \tt Email: wangtao@henu.edu.cn}}
\date{}
\maketitle
\begin{abstract}
A graph is $k$-degenerate if every subgraph $H$ has a vertex $v$ with $d_{H}(v) \leq k$. The class of degenerate graphs plays an important role in the graph coloring theory. Observed that every $k$-degenerate graph is $(k + 1)$-choosable and $(k + 1)$-DP-colorable. Bernshteyn and Lee defined a generalization of $k$-degenerate graphs, which is called \emph{weakly $k$-degenerate}. The weak degeneracy plus one is an upper bound for many graph coloring parameters, such as choice number, DP-chromatic number and DP-paint number. In this paper, we give two sufficient conditions for a plane graph without $4$- and $6$-cycles to be weakly $2$-degenerate, which implies that every such graph is $3$-DP-colorable and near-bipartite, where a graph is near-bipartite if its vertex set can be partitioned into an independent set and an acyclic set. 
\end{abstract}

\section{Introduction}

A graph is {\em $k$-degenerate} if every (induced) subgraph has a vertex of degree at most $k$ in this subgraph. The class of $k$-degenerate graphs plays an important role in the graph coloring theory. Greedy algorithm shows that every $k$-degenerate graph is $(k + 1)$-choosable. For example, every tree is $1$-degenerate, thus it is $2$-choosable. By Euler's formula, every planar graph is a $5$-degenerate graph, and hence it is $6$-choosable. It is well known that not every planar graph is 4-degenerate, but every planar graph is $5$-choosable. 

DP-coloring was introduced in~\cite{MR3758240} by Dvo\v{r}\'{a}k and Postle, it is a generalization of list coloring. Dvo\v{r}\'{a}k and Postle observed that every $k$-degenerate graph is also $(k + 1)$-DP-colorable. 

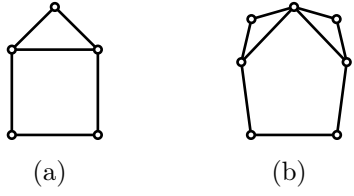
\begin{figure}[htbp]%
\centering
\subcaptionbox{\label{fig:subfig:a-}}
{
\begin{tikzpicture}[scale = 0.8, line width = 1pt]
\coordinate (A) at (45:1);
\coordinate (B) at (135:1);
\coordinate (C) at (225:1);
\coordinate (D) at (-45:1);
\coordinate (H) at (90:1.414);
\draw (A)--(H)--(B)--(C)--(D)--cycle;
\draw (A)--(B);
\node[circle, inner sep = 1, fill = white, draw] () at (A) {};
\node[circle, inner sep = 1, fill = white, draw] () at (B) {};
\node[circle, inner sep = 1, fill = white, draw] () at (C) {};
\node[circle, inner sep = 1, fill = white, draw] () at (D) {};
\node[circle, inner sep = 1, fill = white, draw] () at (H) {};
\end{tikzpicture}}\hspace{1.5cm}
\subcaptionbox{\label{fig:subfig:b-}}
{
\begin{tikzpicture}[scale = 0.8, line width = 1pt]
\coordinate (A) at (30:1);
\coordinate (B) at (150:1);
\coordinate (C) at (225:1);
\coordinate (D) at (-45:1);
\coordinate (H) at (90:1.414);
\coordinate (X) at (60:1.4);
\coordinate (Y) at (120:1.4);
\draw (A)--(X)--(H)--(Y)--(B)--(C)--(D)--cycle;
\draw (A)--(H)--(B);
\node[circle, inner sep = 1, fill = white, draw] () at (A) {};
\node[circle, inner sep = 1, fill = white, draw] () at (B) {};
\node[circle, inner sep = 1, fill = white, draw] () at (C) {};
\node[circle, inner sep = 1, fill = white, draw] () at (D) {};
\node[circle, inner sep = 1, fill = white, draw] () at (H) {};
\node[circle, inner sep = 1, fill = white, draw] () at (X) {};
\node[circle, inner sep = 1, fill = white, draw] () at (Y) {};
\end{tikzpicture}}
\caption{Forbidden configurations in \autoref{3D}.}
\label{NOA}
\end{figure}

By Euler's formula and Handshaking Lemma, every planar graph without $3$-cycles is $3$-degenerate. Wang and Lih~\cite{MR1889505} proved that every planar graph without $5$-cycles is $3$-degenerate. Fijav{\v{z}} \etal~\cite{MR1914478} proved that every planar graph without $6$-cycles is $3$-degenerate. For planar graphs without $4$-cycles, they are not necessarily $3$-degenerate. For example, the line graph of the dodecahedral graph --- icosidodecahedral graph. Liu \etal~\cite{MR4078909} proved that planar graphs without $3$-cycles adjacent to $5$-cycles are $3$-degenerate.  Li and Wang~\cite{Li2019} improved it to the following.
\begin{theorem}[Li and Wang~\cite{Li2019}]\label{3D}
Every planar graph without the configurations in \autoref{NOA} is 3-degenerate. 
\end{theorem}

Sittitrai and Nakprasit~\cite{MR4414782} proved that planar graphs without mutually adjacent 3-, 5-, and 6-cycles are 3-degenerate. Rao and Wang~\cite{MR4230514} proved that every planar graph with neither $4$-, $5$-, $7$-cycles nor triangles at distance less than two is $2$-degenerate. Sittitrai and Nakprasit~\cite{MR4300950} gave some sufficient conditions for planar graphs without $4$- and $5$-cycles to be $2$-degenerate. 
\begin{theorem}[Sittitrai and Nakprasit~\cite{MR4300950}]
\text{}
\begin{enumerate}
\item Every planar graph without $4$-, $5$-, $7$-, and $10$-cycles is $2$-degenerate. 
\item Every planar graph without $4$-, $5$-, $7$-, and $11$-cycles is $2$-degenerate. 
\item Every planar graph without $4$-, $5$-, $8$-, and $10$-cycles is $2$-degenerate. 
\item Every planar graph without $4$-, $5$-, $8$-, and $11$-cycles is $2$-degenerate. 
\end{enumerate}
\end{theorem}

Jumnongnit and Pimpasalee~\cite{MR4268698} gave some sufficient conditions for planar graphs without $4$- and $6$-cycles to be $2$-degenerate. 
\begin{theorem}[Jumnongnit and Pimpasalee~\cite{MR4268698}]
\text{}
\begin{enumerate}
\item Every planar graph without $4$-, $6$-, $8$-, and $10$-cycles is $2$-degenerate. 
\item Every planar graph without $4$-, $6$-, $9$-, and $10$-cycles is $2$-degenerate.
\end{enumerate}
\end{theorem}

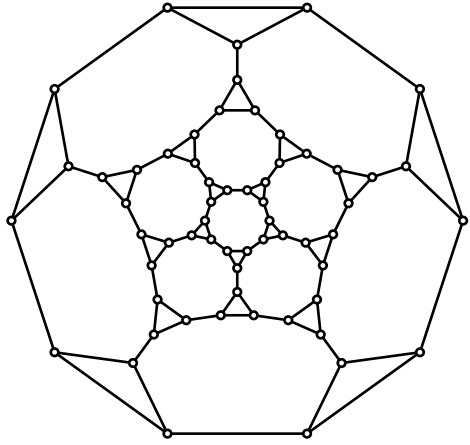
\begin{figure}
\def\s{0.707}
\def\t{0.3}
\def\r{1.2}
\centering
\begin{tikzpicture}[scale = 0.6, line width = 1pt]
\coordinate (O) at (0, 0);
\foreach \ang in {1, 2, 3, 4, 5, 6, 7, 8, 9, 10}
{
\def\pointnamev{v\ang}
\coordinate (\pointnamev) at ($(\ang*360/10:\s)$);
\def\pointnameO{O\ang}
\coordinate (\pointnameO) at ($(\ang*360/10:7*\s)$);
}
\draw (v1)--(v2)--(v3)--(v4)--(v5)--(v6)--(v7)--(v8)--(v9)--(v10)--cycle;
\draw (O1)--(O2)--(O3)--(O4)--(O5)--(O6)--(O7)--(O8)--(O9)--(O10)--cycle;

\coordinate (u1) at ($(v2)!1!60:(v1)$);
\draw (v1)--(u1)--(v2);
\coordinate (u2) at ($(v4)!1!60:(v3)$);
\draw (v3)--(u2)--(v4);
\coordinate (u3) at ($(v6)!1!60:(v5)$);
\draw (v5)--(u3)--(v6);
\coordinate (u4) at ($(v8)!1!60:(v7)$);
\draw (v7)--(u4)--(v8);
\coordinate (u5) at ($(v10)!1!60:(v9)$);
\draw (v9)--(u5)--(v10);

\coordinate (w1) at ($(u1)!1.2!150:(v1)$);
\coordinate (w2) at ($(u2)!1.2!150:(v3)$);
\coordinate (w3) at ($(u3)!1.2!150:(v5)$);
\coordinate (w4) at ($(u4)!1.2!150:(v7)$);
\coordinate (w5) at ($(u5)!1.2!150:(v9)$);
\draw (w1)--(u1);
\draw (w2)--(u2);
\draw (w3)--(u3);
\draw (w4)--(u4);
\draw (w5)--(u5);

\coordinate (x1) at ($(w1)!\r!145:(u1)$);
\coordinate (x2) at ($(w2)!\r!145:(u2)$);
\coordinate (x3) at ($(w3)!\r!145:(u3)$);
\coordinate (x4) at ($(w4)!\r!145:(u4)$);
\coordinate (x5) at ($(w5)!\r!145:(u5)$);

\coordinate (y1) at ($(w1)!\r!-145:(u1)$);
\coordinate (y2) at ($(w2)!\r!-145:(u2)$);
\coordinate (y3) at ($(w3)!\r!-145:(u3)$);
\coordinate (y4) at ($(w4)!\r!-145:(u4)$);
\coordinate (y5) at ($(w5)!\r!-145:(u5)$);

\draw (w1)--(x1)--(y1)--cycle;
\draw (w2)--(x2)--(y2)--cycle;
\draw (w3)--(x3)--(y3)--cycle;
\draw (w4)--(x4)--(y4)--cycle;
\draw (w5)--(x5)--(y5)--cycle;

\coordinate (z1) at (18:4.4*\s);
\coordinate (A1) at (18:5.5*\s);
\coordinate (z11) at ($(z1)!1!150:(A1)$);
\coordinate (z12) at ($(z1)!1!60:(z11)$);
\draw (z1)--(z11)--(z12)--cycle;
\draw (A1)--(z1);
\draw (z11)--(x1);
\draw (z12)--(y5);

\coordinate (z2) at (90:4.4*\s);
\coordinate (A2) at (90:5.5*\s);
\coordinate (z21) at ($(z2)!1!150:(A2)$);
\coordinate (z22) at ($(z2)!1!60:(z21)$);
\draw (z2)--(z21)--(z22)--cycle;
\draw (A2)--(z2);
\draw (z21)--(x2);
\draw (z22)--(y1);

\coordinate (z3) at (162:4.4*\s);
\coordinate (A3) at (162:5.5*\s);
\coordinate (z31) at ($(z3)!1!150:(A3)$);
\coordinate (z32) at ($(z3)!1!60:(z31)$);
\draw (z3)--(z31)--(z32)--cycle;
\draw (A3)--(z3);
\draw (z31)--(x3);
\draw (z32)--(y2);

\coordinate (z4) at (234:4.4*\s);
\coordinate (A4) at (234:5.5*\s);
\coordinate (z41) at ($(z4)!1!150:(A4)$);
\coordinate (z42) at ($(z4)!1!60:(z41)$);
\draw (z4)--(z41)--(z42)--cycle;
\draw (A4)--(z4);
\draw (z41)--(x4);
\draw (z42)--(y3);

\coordinate (z5) at (306:4.4*\s);
\coordinate (A5) at (306:5.5*\s);
\coordinate (z51) at ($(z5)!1!150:(A5)$);
\coordinate (z52) at ($(z5)!1!60:(z51)$);
\draw (z5)--(z51)--(z52)--cycle;
\draw (A5)--(z5);
\draw (z51)--(x5);
\draw (z52)--(y4);

\draw (O10)--(A1)--(O1);
\draw (O2)--(A2)--(O3);
\draw (O4)--(A3)--(O5);
\draw (O6)--(A4)--(O7);
\draw (O8)--(A5)--(O9);

\foreach \ang in {1, 2, 3, 4, 5, 6, 7, 8, 9, 10}
{
\node[circle, inner sep = 1, fill = white, draw] () at (v\ang) {};
\node[circle, inner sep = 1, fill = white, draw] () at (O\ang) {};
}

\foreach \ang in {1, 2, 3, 4, 5}
{
\node[circle, inner sep = 1, fill = white, draw] () at (u\ang) {};
\node[circle, inner sep = 1, fill = white, draw] () at (w\ang) {};
\node[circle, inner sep = 1, fill = white, draw] () at (x\ang) {};
\node[circle, inner sep = 1, fill = white, draw] () at (y\ang) {};
\node[circle, inner sep = 1, fill = white, draw] () at (z\ang) {};
\node[circle, inner sep = 1, fill = white, draw] () at (z\ang1) {};
\node[circle, inner sep = 1, fill = white, draw] () at (z\ang2) {};
\node[circle, inner sep = 1, fill = white, draw] () at (A\ang) {};
}
\end{tikzpicture}
\caption{Truncated dodecahedral graph.}
\label{TD}
\end{figure}


A truncated dodecahedral graph is the graph of vertices and edges of the truncated dodecahedron, see \autoref{TD}. The truncated dodecahedral graph is a $3$-regular planar graph without $4$-, $5$-, $6$-, $7$-, $8$-, and $9$-cycles. It shows that not every planar graphs without $4$-, $5$-, $6$-, $7$-, $8$-, and $9$-cycles is $2$-degenerate. 

Two cycles (or faces) are \emph{normally adjacent} if the intersection is isomorphic to $K_{2}$. In this paper, we prove that every planar graph without $4$-, $6$-, $9$-cycles and a $5$-cycle normally adjacent to a $7$-cycle is almost $2$-degenerate.

\begin{figure}
\centering
\subcaptionbox{A 10-face normally adjacent to a 3-face.\label{10Cap3Fig}}[0.45\linewidth]
{\begin{tikzpicture}[line width = 1pt]
\def\s{0.707}
\foreach \ang in {1, 2, 3, 4, 5, 6, 7, 8, 9, 10}
{
\def\pointname{v\ang}
\coordinate (\pointname) at ($(\ang*360/10+54:\s)$);
}
\draw (v1)--(v2)--(v3)--(v4)--(v5)--(v6)--(v7)--(v8)--(v9)--(v10)--cycle;
\coordinate (S) at ($(v4)!1!60:(v3)$);
\draw (v3)--(S)--(v4);
\node[circle, inner sep = 1, fill, draw] () at (S) {};
\foreach \ang in {1, 2, 3, 4, 5, 6, 7, 8, 9, 10}
{
\node[circle, inner sep = 1, fill, draw] () at (v\ang) {};
}
\end{tikzpicture}}
\subcaptionbox{A bad 10-face adjacent to a special 10-face.\label{BAD10FACE}}[0.45\linewidth]
{\begin{tikzpicture}[line width = 1pt]
\def\s{0.707}
\foreach \ang in {1, 2, 3, 4, 5, 6, 7, 8, 9, 10}
{
\def\pointname{v\ang}
\coordinate (\pointname) at ($(\ang*360/10+54:\s)$);
}
\draw (v1)--(v2)--(v3)--(v4)--(v5)--(v6)--(v7)--(v8)--(v9)--(v10)--cycle;
\foreach \ang in {1, 2, 3, 4, 5, 6, 7, 8, 9, 10}
{
\node[circle, inner sep = 1, fill, draw] () at (v\ang) {};
}

\foreach \ang in {1, 2, 3, 4, 5, 6, 7, 8, 9, 10}
{
\def\pointname{u\ang}
\coordinate (\pointname) at ($(\ang*360/10+54:\s) + (1.9*\s, 0)$);
}
\draw (u1)--(u2)--(u3)--(u4)--(u5)--(u6)--(u7)--(u8)--(u9)--(u10)--cycle;
\foreach \ang in {1, 3, 4, 6, 7, 8, 9, 10}
{
\node[circle, inner sep = 1, fill, draw] () at (u\ang) {};
}
\draw (v10)--(u2);
\draw (v7)--(u5);
\node[rectangle, inner sep = 2, fill, draw] () at (u2) {};
\node[rectangle, inner sep = 2, fill, draw] () at (u5) {};
\end{tikzpicture}}
\caption{Configurations, where a solid dot represents a $3$-vertex, and a rectangle represents a $4$-vertex.}
\label{ReduceC}
\end{figure}
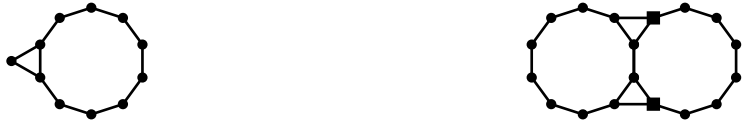

\begin{theorem}\label{469DEGENERATE}
Let $G$ be a plane graph without $4$-, $6$- and $9$-cycles. If there are no $7$-cycles normally adjacent to $5$-cycles, then $\delta(G) \leq 2$ or $G$ contains one of the configurations in \autoref{ReduceC}. 
\end{theorem}
\begin{corollary}[Jumnongnit and Pimpasalee~\cite{MR4268698}]
Every planar graph without $4$-, $6$-, $9$-, and $10$-cycles is $2$-degenerate.
\end{corollary}


A graph is {\em near-bipartite} if its vertex set can be partitioned into an independent set and an acyclic set. Let $\mathscr{G}$ be the class of plane graphs without triangles normally adjacent to $8^{-}$-cycles, without $4$-cycles normally adjacent to $6^{-}$-cycles, and without normally adjacent $5$-cycles. Lu \etal~\cite{MR4422988} proved that every graph in $\mathscr{G}$ is near-bipartite. 

A simple mathematical induction implies that every $2$-degenerate graph is near-bipartite. Although the graphs considered in \autoref{469DEGENERATE} are not $2$-degenerate, the two configurations in \autoref{ReduceC} are reducible configurations for strictly $f$-degenerate transversal (strictly $f$-degenerate transversal was introduced in~\cite{MR4357325}, which is a further generalization of DP-coloring and some vertex-partition problems). As a simple consequence of a result on strictly $f$-degenerate transversal, we have the following result. 
\begin{theorem}
Let $G$ be a plane graph without $4$-, $6$- and $9$-cycles. If there are no $7$-cycles normally adjacent to $5$-cycles, then $G$ is near-bipartite. 
\end{theorem}

In \autoref{sec:3}, we consider another class of planar graphs without $4$-, $6$-, and $8$-cycles. 
\begin{theorem}\label{468DEGENERATE}
Let $G$ be a plane graph without $4$-, $6$-, and $8$-cycles. If there are no $3$-cycles normally adjacent to $9$-cycles, then $\delta(G) \leq 2$ or $G$ contains one of the configurations in \autoref{ReduceC}. 
\end{theorem}

\begin{corollary}[Jumnongnit and Pimpasalee~\cite{MR4268698}]
Every planar graph without $4$-, $6$-, $8$-, and $10$-cycles is $2$-degenerate. 
\end{corollary}

Similar to the above discussions, we can obtain the following result on near-bipartite. 
\begin{theorem}
Let $G$ be a plane graph without $4$-, $6$-, and $8$-cycles. If there are no $3$-cycles normally adjacent to $9$-cycles, then $G$ is near-bipartite. 
\end{theorem}

Motivated by the study of greedy algorithm, Bernshteyn and Lee~\cite{Bernshteyn2021a} defined a new notion of \emph{weak degeneracy}. 
\begin{definition}[\textsf{Delete} Operation]
Let $G$ be a graph and $f : V(G) \longrightarrow \mathbb{Z}$ be a function. For a vertex $u \in V(G)$, the operation \textsf{Delete}$(G, f, u)$ outputs the graph $G' = G - u$ and the function $f': V(G') \longrightarrow \mathbb{Z}$ given by 
\begin{empheq}[left = {f'(v) \coloneqq \empheqlbrace}]{align*}
&f(v) - 1, && \mbox{if $uv \in E(G)$};\\
&f(v), && \mbox{otherwise}.
\end{empheq}
An application of the operation \textsf{Delete} is \emph{legal} if the resulting function $f'$ is nonnegative. 
\end{definition}

\begin{definition}[\textsf{DeleteSave} Operation]
Let $G$ be a graph and $f : V(G) \longrightarrow \mathbb{Z}$ be a function. For a pair of adjacent vertices $u, w \in V(G)$, the operation \textsf{DeleteSave}$(G, f, u, w)$ outputs the graph $G' = G - u$ and the function $f' : V(G') \longrightarrow \mathbb{Z}$ given by
\begin{empheq}[left = {f'(v) \coloneqq \empheqlbrace}]{align*}
&f(v) - 1, &&  \mbox{\quad if $uv \in E(G)$ and $v \neq w$};\\
&f(v), &&  \mbox{\quad otherwise}.
\end{empheq}
An application of the operation \textsf{DeleteSave} is \emph{legal} if $f(u) > f(w)$ and the resulting function $f'$ is nonnegative. 
\end{definition}

A graph $G$ is \emph{weakly $f$-degenerate} if it is possible to remove all vertices from $G$ by a sequence of legal applications of the operations \textsf{Delete} or \textsf{DeleteSave}. A graph is \emph{$f$-degenerate} if it is weakly $f$-degenerate with no \textsf{DeleteSave} operations. Given $d \in \mathbb{N}$, we say that $G$ is \emph{weakly $d$-degenerate} if it is weakly $f$-degenerate with respect to the constant function of value $d$. We say that $G$ is \emph{$d$-degenerate} if it is $f$-degenerate with respect to the constant function of value $d$. The \emph{weak degeneracy} of $G$, denoted by $\mathsf{wd}(G)$, is the minimum integer $d$ such that $G$ is weakly $d$-degenerate. The \emph{degeneracy} of $G$, denoted by $\mathsf{d}(G)$, is the minimum integer $d$ such that $G$ is $d$-degenerate. 

Bernshteyn and Lee~\cite{Bernshteyn2021a} gave the following inequalities. 

\begin{proposition}\label{prop}
For any graph $G$, we always have 
\begin{equation*}
\chi(G) \leq \chi_{\ell}(G) \leq \chi_{\mathsf{DP}}(G) \leq \chi_{\mathsf{DPP}}(G) \leq \mathsf{wd}(G) + 1 \leq \mathsf{d}(G) + 1, 
\end{equation*}
where $\chi_{\mathsf{DP}}(G)$ is the DP-chromatic number of $G$, and $\chi_{\mathsf{DPP}}(G)$ is the DP-paint number of $G$. 
\end{proposition}

From \autoref{prop}, $\mathsf{wd}(G) + 1$ is an upper bound for some graph coloring parameters. So it is interesting to determine the weak degeneracy of a graph. Bernshteyn and Lee~\cite{Bernshteyn2021a} proved that every planar graph is weakly $4$-degenerate, and hence every planar graph is 5-DP-paintable. 

In \autoref{sec:4}, we prove a vertex of degree at most two, and the two configurations in \autoref{ReduceC} are reducible configurations. Then we have the following two results. 

\begin{theorem}
Let $G$ be a plane graph without $4$-, $6$- and $9$-cycles. If there are no $7$-cycles normally adjacent to $5$-cycles, then $G$ is weakly $2$-degenerate. 
\end{theorem}
\begin{theorem}
Let $G$ be a plane graph without $4$-, $6$-, and $8$-cycles. If there are no $3$-cycles normally adjacent to $9$-cycles, then weakly $2$-degenerate. 
\end{theorem}


\begin{corollary}
Let $G$ be a plane graph without $4$-, $6$- and $9$-cycles. If there are no $7$-cycles normally adjacent to $5$-cycles, then $G$ is $3$-DP-colorable, and $3$-DP-paintable. 
\end{corollary}
This strengthens the results that every planar graph without $4$-, $5$-, $6$-, and $9$-cycles is $3$-DP-colorable~\cite{MR3886261}, and every planar graph without $4$-, $6$-, $7$-, and $9$-cycles is $3$-DP-colorable~\cite{MR3969021}. 

\begin{corollary}
Let $G$ be a plane graph without $4$-, $6$-, and $8$-cycles. If there are no $3$-cycles normally adjacent to $9$-cycles, then $G$ is 3-DP-colorable, and $3$-DP-paintable. 
\end{corollary}

We need the following concept in the proofs of \autoref{469DEGENERATE} and \autoref{468DEGENERATE}. In a plane graph, when a face $f$ is normally adjacent to a $3$-face $[uvw]$, where $uv$ is the common edge, we call $w$ a \emph{source} of $f$.  

\section{Proof of \autoref{469DEGENERATE}}
Let $G$ be a minimum counterexample to the statement. Then $G$ is a connected graph with $\delta(G) \geq 3$ and $G$ does not contain any configuration in \autoref{ReduceC}. 

By the structure of the graph, we can easily obtain the following results. 
\begin{lemma}\label{Lem1}
\text{}
\begin{enumerate}[label = (\arabic*)]
    \item A $k$-cycle has no chords, where $k \in \{5, 7\}$.
    \item A $3$-cycle is not adjacent to a $k$-cycle, where $k \in \{3, 5\}$.
    \item The boundaries of two adjacent $5$-faces contain an $8$-cycle.
\end{enumerate}
\end{lemma}

\begin{lemma}\label{Lem2}
\text{}
\begin{enumerate}[label = (\arabic*)]
    \item The boundary of a $6$-face consists of two $3$-cycles. 
    \item The boundary of a $7$-face is a $7$-cycle.
    \item The boundary of an $8$-face consists of an $8$-cycle, or a $3$-cycle and a $5$-cycle, or two $3$-cycles and a cut edge. 
    \item The boundary of a $9$-face consists of three $3$-cycles. 
    \item\label{3F} A $3$-face is not adjacent to a $k$-face, where $k \in \{3, 5, 6, 8, 9\}$.
    \item\label{7F} A $7$-face is adjacent to at most one $3$-face.
    \item\label{No555} A $3$-vertex is incident with at most two $5$-faces. 
    \item\label{5F} A $5$-face is not adjacent to a $k$-face, where $k \in \{6, 7\}$. 
\end{enumerate}
\end{lemma}

A $3$-vertex is \emph{bad} if it is incident with two $5$-faces, \emph{worse} if it is incident with exactly one $5$-face, \emph{worst} if it is incident with a $3$-face. 

\begin{definition}\label{DEF}
Let $f$ be a $10$-face incident with ten $3$-vertices $v_{1}, v_{2}, \dots, v_{10}$. If $f$ is adjacent to two $3$-faces $[v_{1}v_{2}u_{1}]$ and $[v_{3}v_{4}u_{2}]$, then $f$ is called a \emph{bad face}. Since there is no configuration in \autoref{10Cap3Fig}, each of $u_{1}$ and $u_{2}$ is a $4^{+}$-vertex. The edge $v_{2}v_{3}$ is called a \emph{special edge}, $u_{1}v_{2}v_{3}u_{2}$ is called a \emph{special path}, and the face $h$ incident with $u_{1}v_{2}v_{3}u_{2}$ is called a \emph{special face}. 
\end{definition}

By \autoref{Lem2}\ref{3F} and \autoref{Lem2}\ref{7F}, each special face is a $10^{+}$-face incident with at least two $4^{+}$-vertices. 

An initial charge $\mu(x)$ is assigned to each element $x \in V(G) \cup F(G)$ by $\mu(x) = 2d(x) - 6$ for each $x \in V(G)$, and $\mu(x) = d(x) - 6$ for each $x \in F(G)$. By Euler's formula and Handshaking Lemma, the sum of the initial charges is $-12$. Applying the following discharging rules to redistribute the charges, we get a final charge function $\mu'$. We show that $\mu'(x) \geq 0$ for each element $x \in V(G) \cup F(G)$, thus the sum of the final charges is nonnegative, which leads to a contradiction.

\medskip
\textbf{Discharging Rules:}
\begin{enumerate}[label = \textbf{R\arabic*.}, ref = R\arabic*]
    \item\label{R1} Each $3$-face receives $1$ from each incident vertex.
    \item\label{R2} Let $v$ be a $3$-vertex. 
    \begin{enumerate}[label = \textbf{\alph*.}, ref = \alph*]
        \item\label{R2.1} If $v$ is a worst vertex, then it receives $\frac{1}{2}$ from each incident $7^{+}$-face.
        \item\label{R2.2} If $v$ is a worse vertex, then it sends $\frac{1}{2}$ to the incident $5$-face, and receives $\frac{1}{4}$ from each incident $8^{+}$-face. 
        \item\label{R2.3} If $v$ is a bad vertex, then it sends $\frac{1}{8}$ to each incident $5$-face, and receives $\frac{1}{4}$ from the incident $8^{+}$-face. 
    \end{enumerate}
    \item\label{R3} Let $v$ be a $4$-vertex.
    \begin{enumerate}[label = \textbf{\alph*.}, ref = \alph*]
        \item\label{R3.1} If $v$ is not incident with a $3$-face, then $v$ sends $\frac{1}{2}$ to each incident face.
        \item\label{R3.2} If $v$ is incident with a $3$-face and a $5$-face, then it sends $\frac{1}{2}$ to the incident $5$-face, and $\frac{1}{4}$ to each incident $10^{+}$-face. 
        \item\label{R3.3} If $v$ is incident with exactly one $3$-face and no $5$-face, then $v$ sends $\frac{1}{2}$ to each incident face adjacent to the $3$-face. 
    \end{enumerate} 
    \item\label{R4} Each $5$-vertex sends $\frac{2}{3}$ to each incident $5^{+}$-face which is not special. 
    \item\label{R5} Each $6^{+}$-vertex sends $1$ to each incident $5^{+}$-face which is not special. 
    \item\label{R6} Let $\mu^{*}$ denote the charges after applying the rules \ref{R1}--\ref{R5}. Each special face $g$ sends $\frac{\mu^{*}(g)}{t}$ to each adjacent bad face, where $t$ is the number of adjacent bad faces. 
    \item\label{R7} Let $[uvw]$ be a $3$-face, $vw$ be incident with a bad face $f$. 
    \begin{enumerate}[label = \textbf{\alph*.}, ref = \alph*]
        \item\label{R7.1} If $u$ is a $6^{+}$-vertex, and the face $f'$ incident with $uv$ is a special face, then $u$ sends $\frac{1}{2}$ to $f$ via the special face $f'$.  
        \item\label{R7.2} If $u$ is a $5$-vertex, the face $f'$ incident with $uv$ is a special face, and $u$ is the end of two special path in $f'$, then $u$ sends $\frac{1}{3}$ to $f$ via special face $f'$.
        \item\label{R7.3} If $u$ is a $5$-vertex, the face $f'$ incident with $uv$ is a special face, and $u$ is the end of exactly one special path in $f'$, then $u$ sends $\frac{2}{3}$ to $f$ via special face $f'$.
    \end{enumerate}
\end{enumerate}

\begin{lemma}\label{SF}
Let $f$ be a bad face with related labels as in \autoref{DEF}. If both $u_{1}$ and $u_{2}$ are $4$-vertices, then $h$ sends at least $\frac{1}{4}$ to $f$ by \ref{R6}. 
\end{lemma}
\begin{proof}
Assume that $h$ is a $10$-face. If $h$ is incident with exactly two $4^{+}$-vertices, $u_{1}$ and $u_{2}$, then we can easily check that the configuration in \autoref{BAD10FACE} will appear, a contradiction. Hence, $h$ is incident with at least three $4^{+}$-vertices. If $h$ is incident with exactly one special edge, then $\mu^{*}(h) \geq 10 - 6 - 7 \times \frac{1}{2} = \frac{1}{2}$, and $h$ sends at least $\frac{1}{2}$ to $f$. If $h$ is incident with exactly two special edges, then $\mu^{*}(h) \geq 10 - 6 - 7 \times \frac{1}{2} = \frac{1}{2}$, and $h$ sends at least $\frac{1}{4}$ to $f$. If $h$ is incident with exactly three special edges, then $h$ is incident with at least four $4^{+}$-vertices, $\mu^{*}(h) \geq 10 - 6 - 6 \times \frac{1}{2} = 1$, and $h$ sends at least $\frac{1}{3}$ to $f$. Since the degree of $h$ is $10$, it cannot be incident with more than three special edges. 

Assume that $h$ is a $11$-face. If $h$ is incident with exactly one special edge, then $\mu^{*}(h) \geq 11 - 6 - 9 \times \frac{1}{2} = \frac{1}{2}$, and $h$ sends at least $\frac{1}{2}$ to $f$. If $h$ is incident with exactly two special edges, then $h$ is incident with at least three $4^{+}$-vertices, $\mu^{*}(h) \geq 11 - 6 - 8 \times \frac{1}{2} = 1$, and $h$ sends at least $\frac{1}{2}$ to $f$. If $h$ is incident with exactly three special edges, then $h$ is incident with at least four $4^{+}$-vertices, $\mu^{*}(h) \geq 11 - 6 - 7 \times \frac{1}{2} = \frac{3}{2}$, and $h$ sends at least $\frac{1}{2}$ to $f$. Note that $h$ is incident with at most three special edges. 

Assume that $h$ is a $12^{+}$-face. If $h$ is incident with exactly $t$ special edges, then $h$ is incident with at most $d(h) - t$ vertices of degree three, $\mu^{*}(h) \geq d(h) - 6 - (d(h) - t) \times \frac{1}{2} = \frac{d(h) - 12}{2} + \frac{t}{2} \geq \frac{t}{2}$, and $h$ sends at least $\frac{1}{2}$ to $f$.
\end{proof}

Let $v$ be a $3$-vertex. If $v$ is incident with a $3$-face, then the other two incident faces are $7^{+}$-faces by \autoref{Lem2}\ref{3F}, and then $\mu'(v) = -1 + 2 \times \frac{1}{2} = 0$ by \ref{R1} and \ref{R2}\ref{R2.1}. If $v$ is incident with exactly one $5$-face, then the other two incident faces are $8^{+}$-faces by \autoref{Lem2}\ref{5F}, and then $\mu'(v) = -\frac{1}{2} + 2 \times \frac{1}{4} = 0$ by \ref{R2}\ref{R2.2}. If $v$ is incident with exactly two $5$-faces, then the other incident face is an $8^{+}$-face by \autoref{Lem2}\ref{5F}, and then $\mu'(v) = -2 \times \frac{1}{8} + \frac{1}{4} = 0$ \ref{R2}\ref{R2.3}. Since there are no $9$-cycles, $v$ cannot be incident with three $5$-faces by \autoref{Lem2}\ref{No555}. If $v$ is incident with three $6^{+}$-faces, then $\mu'(f) = \mu(f) \geq 0$.

Let $v$ be a $4$-vertex. Then $v$ is incident with at most two $3$-faces. If $v$ is not incident with a $3$-face, then $\mu'(v) = 2 - 4 \times \frac{1}{2} = 0$ by \ref{R3}\ref{R3.1}. If $v$ is incident with two $3$-faces, then $\mu'(v) = 2 - 2 \times 1 = 0$ by \ref{R1}. If $v$ is incident with a $3$-face and a $5$-face, then the other two incident faces are $10^{+}$-faces, and $\mu'(v) = 2 - 1 - \frac{1}{2} - 2 \times \frac{1}{4} = 0$ by \ref{R1} and \ref{R3}\ref{R3.2}. If $v$ is incident with a $3$-face and no $5$-face, then $\mu'(v) = 2 - 1 - 2 \times \frac{1}{2} = 0$ by \ref{R1} and \ref{R3}\ref{R3.3}. 

Let $v$ be a $5$-vertex. Then $v$ is incident with at most two $3$-faces, and $\mu'(v) \geq 4 - 2 \times 1 - 3 \times \frac{2}{3} = 0$ by \ref{R1}, \ref{R4}, \ref{R7}\ref{R7.2} and \ref{R7}\ref{R7.3}. 

Let $v$ be a $6^{+}$-vertex. Then $\mu'(v) \geq \mu(v) - d(v) \times 1 \geq 0$ by \ref{R1}, \ref{R5} and \ref{R7}\ref{R7.1}. 

Assume $f$ is a $d$-face. Since $G$ contains no $4$-cycles, there are no $4$-faces in $G$. 

$\bullet$ $\bm{d = 3}$. Then $\mu'(f) = -3 + 3 \times 1 = 0$ by \ref{R1}.

$\bullet$ $\bm{d = 5}$. If $f$ is incident with at least two $4^{+}$-vertices, then $f$ receives at least $\frac{1}{2}$ from each incident $4^{+}$-vertex by \ref{R3}, \ref{R4}, and \ref{R5}, which implies that $\mu'(f) \geq \mu(f) + 2 \times \frac{1}{2} = 0$. If $f$ is incident with exactly one $4^{+}$-vertex, then $f$ receives at least $\frac{1}{2}$ from the incident $4^{+}$-vertex, and receives at least $\frac{1}{8}$ from each incident $3$-vertex, which implies that $\mu'(f) \geq \mu(f) + \frac{1}{2} + 4 \times \frac{1}{8} = 0$ \ref{R2}. If $f$ is incident with five $3$-vertices, then $f$ is incident with at least one worse vertex by \autoref{Lem2}\ref{No555}, thus $\mu'(f) \geq \mu(f) + \frac{1}{2} + 4 \times \frac{1}{8} = 0$ by \ref{R2}. 

$\bullet$ $\bm{d = 6}$. Then $f$ is not adjacent to a $6^{-}$-face, and $\mu'(f) \geq \mu(f) = 0$.

$\bullet$ $\bm{d = 7}$. Then $f$ is adjacent to at most one $3$-face, and it is not adjacent to a $5$-face. It follows that $\mu'(f) \geq \mu(f) - 2 \times \frac{1}{2} = 0$ by \ref{R2}\ref{R2.1}. 

$\bullet$ $\bm{d = 8}$. Then $f$ is not adjacent to a $3$-face by \autoref{Lem2}\ref{3F}, thus $\mu'(f) \geq \mu(f) - 8 \times \frac{1}{4} = 0$ by \ref{R2}. 

$\bullet$ $\bm{d = 9}$. Then the boundary of $f$ consists of three $3$-cycles, thus it is not adjacent to a $5^{-}$-face. It follows that $\mu'(f) \geq \mu(f) = 3$. 

$\bullet$ $\bm{d = 10}$. Assume $f$ is a special face. Note that $f$ is incident with at least two $4^{+}$-vertices. Then $\mu^{*}(f) \geq 10 - 6 - 8 \times \frac{1}{2} = 0$ by \ref{R2}. By \ref{R6}, its final charge is zero. So we may assume that $f$ is not a special face. Suppose that the boundary is not a $10$-cycle. Then the boundary consists of two $5$-cycles, or a cut edge, a $3$-cycle and a $5$-cycle. In any case, $f$ sends charges to at most eight $3$-vertices, and $\mu'(f) \geq \mu(f) - 8 \times \frac{1}{2} = 0$ by \ref{R2}. So we may assume that the boundary is a $10$-cycle. If $f$ is incident with at most six worst vertices, then $\mu'(f) \geq \mu(f) - 6 \times \frac{1}{2} - 4 \times \frac{1}{4} = 0$ by \ref{R2}. 

{\bfseries Assume that $f$ is incident with exactly seven worst vertices}. Since the integer seven is odd, there is an incident vertex in a $3$-face which is a $4^{+}$-vertex. Then $v$ is incident with at most two worse/bad vertices, and $\mu'(f) \geq \mu(f) - 7 \times \frac{1}{2} - 2 \times \frac{1}{4} = 0$ by \ref{R2}. 

{\bfseries Assume that $f$ is incident with exactly eight worst vertices}. Note that $f$ is not incident with a bad vertex. If $f$ is not incident with a worse vertex, then $\mu'(f) \geq \mu(f) - 8 \times \frac{1}{2} = 0$ by \ref{R2}. 

--- Suppose that $f$ is incident with exactly one worse vertex. Then the other common vertex $w$ on the $5$-face is a $4^{+}$-vertex. If $w$ is a $4$-vertex, then it is not incident with a $3$-face, and $w$ sends $\frac{1}{2}$ to each incident face by \ref{R3}\ref{R3.1}, thus $\mu'(f) \geq \mu(f) - 8 \times \frac{1}{2} - \frac{1}{4} + \frac{1}{2} \geq 0$ by \ref{R2}\ref{R2.1} and \ref{R2}\ref{R2.2}. If $w$ is a $5^{+}$-vertex, then it sends at least $\frac{2}{3}$ to $f$ by \ref{R5} and \ref{R6}, and $\mu'(f) \geq \mu(f) - 8 \times \frac{1}{2} - \frac{1}{4} + \frac{2}{3} \geq 0$ by \ref{R2}\ref{R2.1} and \ref{R2}\ref{R2.2}. 

--- Suppose that $f$ is incident with exactly two worse vertices. Then $f$ is a bad face adjacent to three special faces. Let $v_{i-1}v_{i}$ be incident with a $3$-face $[v_{i-1}v_{i}u_{i/2}]$, where $i \in \{2, 4, 6, 8\}$, and let $v_{9}v_{10}$ be incident with a $5$-face. Since the configuration in \autoref{10Cap3Fig} is forbidden, each of $u_{1}, u_{2}, u_{3}$ and $u_{4}$ is a $4^{+}$-vertex. If $u_{2}$ is a $6^{+}$-vertex, then $u_{2}$ sends at least $2 \times \frac{1}{2} = 1$ to $f$ by \ref{R7}\ref{R7.1}. If $u_{2}$ is a $5$-vertex, then $u_{2}$ is incident with at most two $3$-faces, and $u_{2}$ sends at least $\frac{1}{3} + \frac{2}{3} = 1$ to $f$ by \ref{R7}\ref{R7.2} and \ref{R7}\ref{R7.3}. Hence, if $u_{2}$ is a $5^{+}$-vertex, then $\mu'(f) \geq \mu(f) - 8 \times \frac{1}{2} - 2 \times \frac{1}{4} + 1 > 0$. So we may assume that both $u_{2}$ and $u_{3}$ are $4$-vertices. By \autoref{SF}, the special face incident with $v_{4}v_{5}$ sends at least $\frac{1}{4}$ to $f$. If $u_{1}$ is a $4$-vertex, then the special face incident with $v_{2}v_{3}$ sends at least $\frac{1}{4}$ to $f$ by \autoref{SF}. If $u_{1}$ is a $5^{+}$-vertex, then it sends at least $\frac{1}{3}$ to $f$ by \ref{R7}. Hence, $\mu'(f) \geq 10 - 6 - 8 \times \frac{1}{2} - 2 \times \frac{1}{4} + \min\{2 \times \frac{1}{4}, \frac{1}{4} + \frac{1}{3}\} \geq 0$ whenever $u_{1}$ is a $4^{+}$-vertex. 

{\bfseries Assume that $f$ is incident with exactly nine worst vertices}. Then $f$ is adjacent to five $3$-faces, and it is incident with a $4^{+}$-vertex $w$. Recall that $f$ is not a special face. It follows that $w$ sends at least $\frac{1}{2}$ to $f$ by \ref{R3}\ref{R3.3}, \ref{R4} and \ref{R5}, thus $\mu'(f) \geq \mu(f) - 9 \times \frac{1}{2} + \frac{1}{2} = 0$ by \ref{R2}. 

{\bfseries Assume that $f$ is incident with ten worst vertices}. Then the five sources are $4^{+}$-vertices. If one of the sources is a $5^{+}$-vertex, then it sends at least $1$ to $f$ by \ref{R7}, and $\mu'(f) \geq 10 - 6 - 10 \times \frac{1}{2} + 1 = 0$. So we may assume that the five sources are $4$-vertices. By \autoref{SF}, each of the adjacent special face sends at least $\frac{1}{4}$ to $f$, thus $\mu'(f) \geq 10 - 6 - 10 \times \frac{1}{2} + 5 \times \frac{1}{4} > 0$. 

$\bullet$ $\bm{d = 11}$. Similar to the case that $d = 10$, if $f$ is a special $11$-face, then its final charge is zero by \ref{R6}. So we may assume that $f$ is not a special face. If $f$ is incident with at most nine worst vertices, then $\mu'(f) \geq 11 - 6 - 9 \times \frac{1}{2} - 2 \times \frac{1}{4} = 0$ by \ref{R2}. Since the integer $11$ is odd, $f$ is incident with at most ten worst vertices. So we may assume that $f$ is incident with exactly ten worst vertices. Then the remaining vertex is a $4^{+}$-vertex, or a $3$-vertex incident with three $7^{+}$-faces. It follows that $\mu'(f) \geq 11 - 6 - 10 \times \frac{1}{2} = 0$ by \ref{R2}. 

$\bullet$ $\bm{d \geq 12}$. Similar to the case $d = 10$, if $f$ is a special face, then its final charge is zero by \ref{R6}. So we may assume that $f$ is not a special face. Then $\mu'(f) \geq d - 6 - d \times \frac{1}{2} \geq 0$. 

\section{Proof of \autoref{468DEGENERATE}}\label{sec:3}
Let $G$ be a minimum counterexample to the statement. Then $G$ is a connected graph with $\delta(G) \geq 3$ and $G$ does not contain any configuration in \autoref{ReduceC}. 

By the structure of the graph, we can easily obtain the following properties. 
\begin{enumerate}[label = \textbf{(\arabic*)}, ref = (\arabic*)]
\item The minimum degree is at least three. 
\item\label{4FACE} There are no $4$-faces.
\item\label{NoChord} A $k$-cycle has no chords, where $k \in \{5, 7, 9\}$. 
\item\label{3V5} A $3$-cycle is not adjacent to a $k$-cycle, where $k \in \{3, 5, 7\}$. 
\item There are no adjacent $5$-cycles. 
\item The boundary of a $6$-face consists of two $3$-cycles. 
\item\label{7FACE} The boundary of a $7$-face is a $7$-cycle. 
\item\label{8FACE} The boundary of an $8$-face consists of a $3$-cycle and a $5$-cycle, or two $3$-cycles and a cut edge. 
\end{enumerate}
\begin{enumerate}[label = \textbf{(\arabic*)}, ref = (\arabic*), resume]
\item\label{3V8-} Every $3$-face is adjacent to three $10^{+}$-faces.
\end{enumerate}
\begin{proof}[Proof of \ref{3V8-}]
The boundaries of two adjacent $3$-faces form a $4$-cycle with a chord, this contradicts the hypothesis. So there are no two adjacent $3$-faces. By \ref{4FACE}, there are no $4$-faces. By \ref{3V5}, there are no $3$-faces adjacent to $5$-faces. Note that the boundary of every $6$-face consists of two $3$-cycles, thus there are no $3$-faces adjacent to $6$-faces. By \ref{NoChord} and \ref{7FACE}, if a $3$-face is adjacent to a $7$-face, then there is a $3$-cycle adjacent to an induced $7$-cycle, this contradicts \ref{3V5}. Suppose that a $3$-face $f$ is adjacent to an $8$-face $g$. By \ref{8FACE}, every edge on the boundary of $g$ is a cut edge, or in a $3$- or $5$-cycle, so it cannot be included in another $3$-cycle which is the boundary of $f$. It is easy to check that the boundary of every $9$-face consists of three triangles or a $9$-cycle. If the boundary of a $9$-face $f$ is a $9$-cycle, then it is an induced $9$-cycle by \ref{NoChord}, and $f$ cannot be adjacent to a $3$-face. If the boundary of a $9$-face $f$ consists of three triangles, then no edge on the boundary of $f$ can be included in other $3$-cycles. Therefore, every $3$-face cannot be adjacent to any $9^{-}$-face. 
\end{proof}

\begin{enumerate}[label = \textbf{(\arabic*)}, ref = (\arabic*), resume]
\item\label{5V6-} Every face adjacent to a $5$-face is a $7^{+}$-face. 
\end{enumerate}
\begin{proof}[Proof of \ref{5V6-}]
Let $f = [v_{1}v_{2}v_{3}v_{4}v_{5}]$ be a $5$-face. Note that $f$ cannot be adjacent to $4^{-}$-faces. Suppose that $g = [v_{1}u_{2}u_{3}u_{4}v_{5}]$ is a $5$-face adjacent to $f$. Since the minimum degree is at least three and every $5$-cycle has no chord, we have that $v_{1}, v_{2}, v_{3}, v_{4}, v_{5}, u_{2}, u_{4}$ are seven distinct vertices. Since $f$ and $g$ are symmetric, we have that $u_{3}$ does not belong to $\{v_{2}, v_{4}\}$. It is observed that $u_{3} \neq v_{3}$, for otherwise $[v_{3}v_{4}v_{5}u_{4}]$ is a $4$-cycle. Hence, the boundaries of $f$ and $g$ are normally adjacent, but they form an $8$-cycle with a chord, a contradiction. 

Suppose that $f$ is adjacent to a $6$-face. Since the boundary of a $6$-face consists of two $3$-cycles, we have that the boundary of $f$ is adjacent to a $3$-cycle, this contradicts \ref{3V5}. Therefore, every face adjacent to $f$ is a $7^{+}$-face. 
\end{proof}

A $5$-face is \emph{nice} if it is incident with five $3$-vertices. A $3$-vertex is \emph{poor} if it is incident with a nice $5$-face, \emph{worst} if it is incident with a $3$-face. The definitions of \emph{bad face}, \emph{special edge}, \emph{special path} and \emph{special face} are the same with that in \autoref{DEF}. 

As in the proof of \autoref{469DEGENERATE}, we still use the discharging method with initial charge $\mu(x) = 2d(x) - 6$ for each $x \in V(G)$, and $\mu(x) = d(x) - 6$ for each $x \in F(G)$. Let $\mu'(x)$ to denote the final charge after the following nine rules are applied. 

\medskip
\textbf{Discharging Rules:}
\begin{enumerate}[label = \textbf{R\arabic*.}, ref = R\arabic*]
    \item\label{468R1} Each $3$-face receives $1$ from each incident vertex.
    \item\label{468R2} Each $5$-face receives $1$ from each incident $4^{+}$-vertex. 
    \item\label{468R3} Each nice $5$-face receives $\frac{1}{5}$ from each incident $3$-vertex.
    \item\label{468R4} Let $v$ be a $3$-vertex. 
    \begin{enumerate}[label = \textbf{\alph*.}, ref = \alph*]
        \item\label{468R4.1} If $v$ is a worst vertex, then it receives $\frac{1}{2}$ from each incident $10^{+}$-face.
        \item\label{468R4.2} If $v$ is a poor vertex, then it receives $\frac{1}{10}$ from each incident $7^{+}$-face. 
    \end{enumerate}
    \item\label{468R5} Let $v$ be a $4$-vertex.
    \begin{enumerate}[label = \textbf{\alph*.}, ref = \alph*]
        \item\label{468R5.1} If $v$ is not incident with a $5^{-}$-face, then $v$ sends $\frac{1}{2}$ to each incident face.
        \item\label{468R5.2} If $v$ is incident with exactly one $5^{-}$-face, then $v$ sends $\frac{1}{2}$ to each incident face adjacent to the $5^{-}$-face. 
    \end{enumerate} 
    \item\label{468R6} Each $5$-vertex sends $\frac{2}{3}$ to each incident non-special $6^{+}$-face. 
    \item\label{468R7} Each $6^{+}$-vertex sends $1$ to each incident non-special $6^{+}$-face. 
    \item\label{468R8} Let $\mu^{*}$ denote the charges after applying the rules \ref{468R1}--\ref{468R7}. Each special face $g$ sends $\frac{\mu^{*}(g)}{t}$ to each adjacent bad face, where $t$ is the number of adjacent bad faces. 
    \item\label{468R9} Let $[uvw]$ be a $3$-face, $vw$ be incident with a bad face $f$. 
    \begin{enumerate}[label = \textbf{\alph*.}, ref = \alph*]
        \item\label{468R9.1} If $u$ is a $6^{+}$-vertex, and the face $f'$ incident with $uv$ is a special face, then $u$ sends $\frac{1}{2}$ to $f$ via the special face $f'$.  
        \item\label{468R9.2} If $u$ is a $5$-vertex, the face $f'$ incident with $uv$ is a special face, and $u$ is the end of two special path in $f'$, then $u$ sends $\frac{1}{3}$ to $f$ via special face $f'$.
        \item\label{468R9.3} If $u$ is a $5$-vertex, the face $f'$ incident with $uv$ is a special face, and $u$ is the end of exactly one special path in $f'$, then $u$ sends $\frac{2}{3}$ to $f$ via special face $f'$.
    \end{enumerate}
\end{enumerate}

\begin{lemma}\label{468SF}
Let $f$ be a bad face with related labels as in \autoref{DEF}. If both $u_{1}$ and $u_{2}$ are $4$-vertices, then $h$ sends at least $\frac{1}{4}$ to $f$ by \ref{468R8}. 
\end{lemma}
\begin{proof}
Assume that $h$ is a $10$-face. If $h$ is incident with exactly two $4^{+}$-vertices, $u_{1}$ and $u_{2}$, then we can easily check that the configuration in \autoref{BAD10FACE} will appear, a contradiction. Hence, $h$ is incident with at least three $4^{+}$-vertices. If $h$ is incident with exactly one special edge, then $\mu^{*}(h) \geq 10 - 6 - 7 \times \frac{1}{2} = \frac{1}{2}$, and $h$ sends at least $\frac{1}{2}$ to $f$. If $h$ is incident with exactly two special edges, then $\mu^{*}(h) \geq 10 - 6 - 7 \times \frac{1}{2} = \frac{1}{2}$, and $h$ sends at least $\frac{1}{4}$ to $f$. If $h$ is incident with exactly three special edges, then $h$ is incident with at least four $4^{+}$-vertices, $\mu^{*}(h) \geq 10 - 6 - 6 \times \frac{1}{2} = 1$, and $h$ sends at least $\frac{1}{3}$ to $f$. Since the degree of $h$ is $10$, it cannot be incident with more than three special edges. 

Assume that $h$ is a $11$-face. If $h$ is incident with exactly one special edge, then $\mu^{*}(h) \geq 11 - 6 - 9 \times \frac{1}{2} = \frac{1}{2}$, and $h$ sends at least $\frac{1}{2}$ to $f$. If $h$ is incident with exactly two special edges, then $h$ is incident with at least three $4^{+}$-vertices, $\mu^{*}(h) \geq 11 - 6 - 8 \times \frac{1}{2} = 1$, and $h$ sends at least $\frac{1}{2}$ to $f$. If $h$ is incident with exactly three special edges, then $h$ is incident with at least four $4^{+}$-vertices, $\mu^{*}(h) \geq 11 - 6 - 7 \times \frac{1}{2} = \frac{3}{2}$, and $h$ sends at least $\frac{1}{2}$ to $f$. Note that $h$ is incident with at most three special edges. 

Assume that $h$ is a $12^{+}$-face. If $h$ is incident with exactly $t$ special edges, then $h$ is incident with at most $d(h) - t$ vertices of degree three, $\mu^{*}(h) \geq d(h) - 6 - (d(h) - t) \times \frac{1}{2} = \frac{d(h) - 12}{2} + \frac{t}{2} \geq \frac{t}{2}$, and $h$ sends at least $\frac{1}{2}$ to $f$.
\end{proof}

Let $v$ be a $3$-vertex. If $v$ is incident with a $3$-face, then the other two incident faces are $10^{+}$-faces by \ref{3V8-}, and then $\mu'(v) = -1 + 2 \times \frac{1}{2} = 0$ by \ref{468R1} and \ref{468R4}\ref{468R4.1}. If $v$ is a poor vertex, then the other two incident faces are $7^{+}$-faces by \ref{5V6-}, and then $\mu'(v) = -\frac{1}{5} + 2 \times \frac{1}{10} = 0$ by \ref{468R3} and \ref{468R4}\ref{468R4.2}. If $v$ is not a worst or poor vertex, then it is not involved in the discharging rules, $\mu'(v) = \mu(v) = 0$. 

Let $v$ be a $4$-vertex. Then $v$ is incident with at most two $5^{-}$-faces. If $v$ is not incident with a $5^{-}$-face, then $\mu'(v) = 2 - 4 \times \frac{1}{2} = 0$ by \ref{468R5}\ref{468R5.1}. If $v$ is incident with exactly one $5^{-}$-face, then $\mu'(v) = 2 - 1 - 2 \times \frac{1}{2} = 0$ by \ref{468R1}, \ref{468R2} and \ref{468R5}\ref{468R5.2}. If $v$ is incident with two $5^{-}$-faces, then $\mu'(v) = 2 - 2 \times 1 = 0$ by \ref{468R1} and \ref{468R2}. 

Let $v$ be a $5$-vertex. Then $v$ is incident with at most two $5^{-}$-faces. By \ref{468R1} and \ref{468R2}, $v$ sends $1$ to each incident $5^{-}$-face. By \ref{468R6}, \ref{468R9}\ref{468R9.2} and \ref{468R9}\ref{468R9.3}, $v$ sends at most $\frac{2}{3}$ to/via each incident $6^{+}$-face. It follows that $\mu'(v) \geq 4 - 2 \times 1 - 3 \times \frac{2}{3} = 0$. 

Let $v$ be a $6^{+}$-vertex. Then $\mu'(v) \geq \mu(v) - d(v) \times 1 \geq 0$ by \ref{468R1}, \ref{468R2}, \ref{468R7} and \ref{468R9}\ref{468R9.1}.

Assume $f$ is a $d$-face. Since $G$ contains no $4$-cycles, there are no $4$-faces in $G$. 

$\bullet$ $\bm{d = 3}$. Then $\mu'(f) = -3 + 3 \times 1 = 0$ by \ref{468R1}.

$\bullet$ $\bm{d = 5}$. If $f$ is incident with a $4^{+}$-vertex, then $\mu'(f) \geq \mu(f) + 1 = 0$ by \ref{468R2}. If $f$ is not incident with a $4^{+}$-vertex, then it is incident with five $3$-vertices, and $\mu'(f) = \mu(f) + 5 \times \frac{1}{5} = 0$ by \ref{468R3}.  

$\bullet$ $\bm{d = 6}$. By \ref{3V8-} and \ref{5V6-}, $f$ is not adjacent to a $5^{-}$-face, and $\mu'(f) \geq \mu(f) = 0$.

$\bullet$ $\bm{d \in \{7, 8, 9\}}$. By \ref{3V8-}, $f$ is not adjacent to a $3$-face. It follows that $\mu'(f) \geq d - 6 - d \times \frac{1}{10} > 0$ by \ref{468R4}\ref{468R4.2}. 

$\bullet$ $\bm{d = 10}$. Assume that $f$ is a special face. Note that $f$ is incident with at least two $4^{+}$-vertices. Then $\mu^{*}(f) \geq 10 - 6 - 8 \times \frac{1}{2} = 0$ by \ref{468R4}. By \ref{468R8}, its final charge is zero.  So we may assume that $f$ is not a special face. If $f$ is incident with at most seven worst vertices, then $\mu'(f) \geq \mu(f) - 7 \times \frac{1}{2} - 3 \times \frac{1}{10} > 0$ by \ref{468R4}. 

{\bfseries Assume that $f$ is incident with exactly eight worst vertices}. We may further assume that $f$ is adjacent to a nice $5$-face, otherwise $f$ is not adjacent to a nice $5$-face, and $\mu'(f) \geq \mu(f) - 8 \times \frac{1}{2} = 0$ by \ref{468R4}.  Since the configuration in \autoref{10Cap3Fig} is forbidden, all the four sources are $4^{+}$-vertices. If one of the sources is a $5^{+}$-vertex, then it sends at least $\frac{1}{3}$ to $f$ by \ref{468R9}, thus $\mu'(f) \geq \mu(f) - 8 \times \frac{1}{2} - 2 \times \frac{1}{10} + \frac{1}{3} > 0$ by \ref{468R4}. If all the four sources are $4$-vertices, then $f$ is adjacent to three special faces, which implies that $\mu'(f) \geq \mu(f) - 8 \times \frac{1}{2} - 2 \times \frac{1}{10} + 3 \times \frac{1}{4} > 0$ by \ref{468R4} and \autoref{468SF}. 

{\bfseries Assume that $f$ is incident with exactly nine worst vertices}. Then $f$ is adjacent to five $3$-faces, and it is incident with a $4^{+}$-vertex $w$. Recall that $f$ is not a special face. It follows that $w$ sends at least $\frac{1}{2}$ to $f$ by \ref{468R5}\ref{468R5.2}, \ref{468R6} and \ref{468R7}, thus $\mu'(f) \geq \mu(f) - 9 \times \frac{1}{2} + \frac{1}{2} = 0$ by \ref{468R4}. 

{\bfseries Assume that $f$ is incident with ten worst vertices}. Then the five sources are $4^{+}$-vertices. If one of the sources is a $5^{+}$-vertex, then it sends at least $1$ to $f$ by \ref{468R9}, and $\mu'(f) \geq 10 - 6 - 10 \times \frac{1}{2} + 1 = 0$. So we may assume that the five sources are $4$-vertices. By \autoref{468SF}, each of the adjacent special face sends at least $\frac{1}{4}$ to $f$, thus $\mu'(f) \geq 10 - 6 - 10 \times \frac{1}{2} + 5 \times \frac{1}{4} \geq 0$. 

$\bullet$ $\bm{d = 11}$. Similar to the case that $d = 10$, if $f$ is a special face, then its final charge is zero by \ref{468R8}. So we may assume that $f$ is not a special face. If $f$ is incident with at most nine worst vertices, then $\mu'(f) \geq 11 - 6 - 9 \times \frac{1}{2} - 2 \times \frac{1}{10} > 0$ by \ref{468R4}. Since the integer $11$ is odd, $f$ is incident with at most ten worst vertices. So we may assume that $f$ is incident with exactly ten worst vertices. Then the remaining vertex is a $4^{+}$-vertex, or a $3$-vertex incident with three $7^{+}$-faces. It follows that $\mu'(f) \geq 11 - 6 - 10 \times \frac{1}{2} = 0$ by \ref{468R4}. 

$\bullet$ $\bm{d \geq 12}$. Similar to the case $d = 10$, if $f$ is a special face, then its final charge is zero by \ref{468R8}. So we may assume that $f$ is not a special face. Then $\mu'(f) \geq d - 6 - d \times \frac{1}{2} \geq 0$. 

\section{Reducible configurations for weakly 2-degenerate}\label{sec:4}
In this section, we will prove three reducible configurations. We say that $G$ is a \emph{minimal} graph of weak degeneracy $d$ if $\mathsf{wd}(G) = d$ and $\mathsf{wd}(H) < d$ for every proper subgraph $H$ of $G$. A connected graph is a \emph{GDP-tree} if every block is either a cycle or a complete graph. Bernshteyn and Lee~\cite{Bernshteyn2021a} proved the following Gallai type result.

\begin{theorem}\label{Gallai}
Let $G$ be a minimal graph of weak degeneracy $d \geq 3$. 
\begin{enumerate}[label = (\roman*)]
\item\label{G1} The minimum degree of $G$ is at least $d$. 
\item\label{G2} Let $U \subseteq \{u \in V(G) \mid d_{G}(u) = d\}$. Then every component of $G[U]$ is a GDP-tree.
\end{enumerate}
\end{theorem}

Assume $G$ is a graph which is not weakly $2$-degenerate but every proper subgraph is weakly $2$-degenerate. Then $G$ is a connected graph with $\mathsf{wd}(G) = 3$. Hence, $G$ is a minimal graph of weak degeneracy 3. By \autoref{Gallai}\ref{G1}, a vertex of degree at most two is reducible. 

Let $G'$ be the subgraph induced by the vertices represented in \autoref{10Cap3Fig}. The graph $G'$ contains the graph represented in \autoref{10Cap3Fig} as a spanning subgraph, so $G'$ is a 2-connected graph which is neither a cycle nor a complete graph. By \autoref{Gallai}\ref{G2}, the graph represented in \autoref{10Cap3Fig} cannot be a subgraph of a minimal graph of weak degeneracy 3.

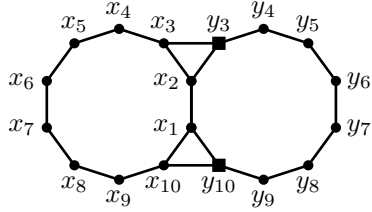
\begin{figure}
\centering
\begin{tikzpicture}[line width = 1pt]
\def\s{1}
\foreach \ang in {1, 2, 3, 4, 5, 6, 7, 8, 9, 10}
{
\def\pointname{v\ang}
\coordinate (\pointname) at ($(\ang*360/10+54:\s)$);
}
\draw (v1)node[above]{\small$x_{4}$}--(v2)node[above]{\small$x_{5}$}--(v3)node[left]{\small$x_{6}$}--(v4)node[left]{\small$x_{7}$}--(v5)node[below]{\small$x_{8}$}--(v6)node[below]{\small$x_{9}$}--(v7)node[below]{\small$x_{10}$}--(v8)node[left]{\small$x_{1}$}--(v9)node[left]{\small$x_{2}$}--(v10)node[above]{\small$x_{3}$}--cycle;
\foreach \ang in {1, 2, 3, 4, 5, 6, 7, 8, 9, 10}
{
\node[circle, inner sep = 1, fill, draw] () at (v\ang) {};
}

\foreach \ang in {1, 2, 3, 4, 5, 6, 7, 8, 9, 10}
{
\def\pointname{u\ang}
\coordinate (\pointname) at ($(\ang*360/10+54:\s) + (1.9*\s, 0)$);
}
\draw (u1)node[above]{\small$y_{4}$}--(u2)node[above]{\small$y_{3}$}--(u3)--(u4)--(u5)node[below]{\small$y_{10}$}--(u6)node[below]{\small$y_{9}$}--(u7)node[below]{\small$y_{8}$}--(u8)node[right]{\small$y_{7}$}--(u9)node[right]{\small$y_{6}$}--(u10)node[above]{\small$y_{5}$}--cycle;
\foreach \ang in {1, 3, 4, 6, 7, 8, 9, 10}
{
\node[circle, inner sep = 1, fill, draw] () at (u\ang) {};
}
\draw (v10)--(u2);
\draw (v7)--(u5);
\node[rectangle, inner sep = 2, fill, draw] () at (u2) {};
\node[rectangle, inner sep = 2, fill, draw] () at (u5) {};
\end{tikzpicture}
\caption{A bad 10-face adjacent to a special 10-face.}
\label{CF}
\end{figure}

Let $W$ be the set of vertices represented in \autoref{CF}. Since $G$ contains no $4$- and $6$-cycles, the edges incident with $x_{4}$ in $G[W]$ are $x_{4}x_{3}, x_{4}x_{5}$ and possible edges $x_{4}x_{6}, x_{4}x_{8}, x_{4}y_{7}$ or $x_{4}y_{8}$. As $x_{4}$ is a 3-vertex, at most one of $x_{4}x_{6}, x_{4}x_{8}, x_{4}y_{7}$ and $x_{4}y_{8}$ can appear in $G[W]$. If $x_{4}x_{8} \in E(G)$, then there exist an $8$-cycle $x_{1}x_{2}x_{3}x_{4}x_{8}x_{9}x_{10}y_{10}x_{1}$ and a $9$-cycle $x_{1}x_{2}y_{3}x_{3}x_{4}x_{8}x_{9}x_{10}y_{10}x_{1}$. If $x_{4}y_{7} \in E(G)$, then there exist an $8$-cycle $x_{4}y_{7}y_{8}y_{9}y_{10}x_{1}x_{2}x_{3}x_{4}$ and a $9$-cycle $x_{4}y_{7}y_{8}y_{9}y_{10}x_{10}x_{1}x_{2}x_{3}x_{4}$. If $x_{4}y_{8} \in E(G)$, then there exist an $8$-cycle $x_{4}y_{8}y_{9}y_{10}x_{10}x_{1}x_{2}x_{3}x_{4}$ and a $9$-cycle $x_{4}y_{8}y_{9}y_{10}x_{10}x_{1}x_{2}y_{3}x_{3}x_{4}$. Note that either $8$-cycles or $9$-cycles are forbidden in \autoref{469DEGENERATE} and \autoref{468DEGENERATE}, then $x_{4}x_{8}, x_{4}y_{7}, x_{4}y_{8} \notin E(G)$. Hence, $x_{4}$ may have only one possible edge $x_{4}x_{6}$. By the symmetry of $x_{4}$ and $x_{9}$, the vertex $x_{9}$ has only one possible edge $x_{9}x_{7}$. When $x_{4}x_{6}, x_{9}x_{7} \in E(G)$, there exist an $8$-cycle $x_{1}x_{2}x_{3}x_{4}x_{6}x_{7}x_{9}x_{10}x_{1}$ and a $9$-cycle $x_{1}x_{2}x_{3}x_{4}x_{6}x_{7}x_{8}x_{9}x_{10}x_{1}$, a contradiction. Assume, without loss of generality, that $x_{4}x_{6} \notin E(G)$. 

By minimality, we can remove all vertices from $G - W$ by a sequence of legal applications of the operations \textsf{Delete} and \textsf{DeleteSave}, and then we remove the vertices from $W$ with the ordering 
\[
x_{3}, y_{3}, y_{4}, y_{5}, y_{6}, y_{7}, y_{8}, y_{9}, y_{10}, x_{2}, x_{1}, x_{10}, x_{9}, x_{8}, x_{7}, x_{6}, x_{5}, x_{4}
\]
by a sequence of legal applications of the operations \textsf{Delete} except the first operation \textsf{DeleteSave}$(G, \cdot\,, x_{3}, x_{4})$. Hence, the graph represented in \autoref{BAD10FACE} cannot be a subgraph of a counterexample to \autoref{469DEGENERATE} or \autoref{468DEGENERATE}. 

Therefore, the graphs considered in \autoref{469DEGENERATE} and \autoref{468DEGENERATE} are weakly $2$-degenerate.

\end{document}